\documentclass[12pt, a4paper]{scrartcl}

\usepackage{times} % times is used to avoid bitmap fonts in PDF
\usepackage{amssymb}
\usepackage{amsmath}
\usepackage{amscd}
 
\newcommand{\qed}{\mbox{}\hspace*{\fill}\nolinebreak\mbox{$\rule{0.6em}{0.6em}$}} %%to end your proof write $\qed$.
\begin{document}

\begin{center}

{\large The Unitary Group In Its Strong Topology}\\
\vspace{0,2cm}
Martin Schottenloher \\
Mathematisches Institut
LMU M\"unchen\\
Theresienstr. 39, 80333 M\"unchen\\
\texttt{schotten@math.lmu.de}, +49 89 21804435 

\end{center}
\vspace{1,0cm}

\begin{abstract} 
\textbf{Abstract.} The unitary group $\mathrm U(\mathcal H)$ on an infinite dimensional complex Hilbert space $\mathcal H$ in its strong topology is a topological group and has some further nice properties, e.g. it is metrizable and contractible if $\mathcal H$ is separable. As an application Hilbert bundles are classified by homotopy. 
\end{abstract}

\begin{abstract} 
\end{abstract}

\begin{abstract} \textit{2000 Mathematics Subject Classification:} 51 F 25, 22 A 10, 22 E 65, 55 R 35, 57 T 20
\end{abstract}
	
\vspace{1cm}
	
\noindent It is easy to show and well-known that the unitary group $\mathrm U(\mathcal H)$ 
-- the group of all unitary operators $\mathcal H \to \mathcal H$ on a complex Hilbert space $\mathcal H$ -- 
is a topological group with respect to the norm topology on $\mathrm U(\mathcal H)$.
However, for many purposes in mathematics the norm topology is too strong. For example, for a compact topological group $G$ with Haar measure $\mu$ the left regular representation on $\mathcal H = L_2(G,\mu)$
	$$L : G \to \mathrm U(\mathcal H), g\mapsto L_g: \mathcal H \to \mathcal H, L_gf(x)= f(g^{-1}x)$$ 
is continuous for the strong topology on $\mathrm U(\mathcal H)$, but $L$ is not continuous when $\mathrm U(\mathcal H)$ is equipped with the norm topology, except for finite $G$. This fact makes the norm topology on $\mathrm U(\mathcal H)$ useless in representation theory and its applications as well as in many areas of physics or topology. The continuity property which is mostly used in case of a topological space $W$ and a general Hilbert space $\mathcal H$ and which seems to be more natural is the continuity of a left action on $\mathcal H$
	$$\Phi: W\times \mathcal H \to \mathcal H,$$
in particular, in case of a left action of a topological group $G$ on $\mathcal H$: Note that the above left regular representation is continuous as a map: $L : G\times \mathcal H \to \mathcal H$. \\

Whenever $\Phi$ is a unitary action (i.e. $\widehat \Phi(w): f\mapsto \Phi(w,f)$ is a unitary operator $\widehat \Phi(w) \in \mathrm U(\mathcal H)$ for all $w\in W$) the continuity of $\Phi$ is equivalent to the continuity of the induced map
	$$\widehat \Phi: W \to \mathrm U(\mathcal H)$$
with respect to the strong topology on $\mathrm U(\mathcal H)$. In fact, if the action  $\Phi$ is continuous then $\widehat \Phi$ is strongly continuous by definition of the strong topology. The converse holds since $\mathrm U(\mathcal H)$ is a uniformly bounded set of operators. The corresponding statement for the general linear group $\mathrm{GL} (\mathcal H)$ of bounded invertible operators holds for the compact open topology on $\mathrm{GL} (\mathcal H)$ instead of the strong topology. On $\mathrm U(\mathcal H)$ the two topologies coincide, see proposition 2 below.\\

We come back to this property in a broader context at the end of this short paper where we explain the significance of the fact that $\mathrm{U}(\mathcal H)$ is a topological group for the classification of Hilbert bundles over paracompact spaces $X$.\\
%Another example which shows the advantage of the strong topology in comparison to the norm topology on $\mathrm U(\mathcal H)$ or the projective unitary group $\mathrm{PU}(\mathcal H)$ is the following. If $g: Y \to X$ is a fiber bundle with compact but infinite fibers $Y_x = g^{-1}(x), x\in Y,$ the induced bundle $F\to X$ whose fibers $F_x$ are the $L_2$-spaces of half-densities on $Y_x$ (resp. the corresponding projective spaces) is a bundle of infinite dimensional Hilbert spaces (resp. of its projectivations) on which $\mathrm{U}(\mathcal H)$ (resp. $\mathrm{PU}(\mathcal H)$) acts naturally. The action is  continuous in the strong topology while it is not continuous in the norm topology. \\

\textbf{Proposition 1:} \textit{$\mathrm U(\mathcal H)$ is a topological group with respect to the strong topology.}\\

\textbf{Proof.} Indeed, the composition $(S,T) \mapsto ST$ is continuous:  
Given $(S_0,T_0) \in \mathrm U(\mathcal H)$ let $\mathcal V$ be a neighbourhood of $S_0T_0$ of the form 
$\mathcal V = \left\{R \in \mathcal H : \left\|(R - S_0T_0)f\right\| < \varepsilon\right\}$ where
$f\in \mathcal H$ and $\varepsilon > 0$. Now, $\mathcal W := \left\{(S,T) \in \mathcal H: \left\|(S - S_0)T_0f\right\| < \frac{1}{2}\varepsilon, \left\|(T - T_0)f\right\| < \frac{1}{2}\varepsilon \right\}$ is a neighbourhood of $(S_0,T_0)\in \mathcal H$ and for $(S,T)\in \mathcal W$ we have: 
$$\left\|(ST - S_0T_0)f\right\| \leq \left\|S(T - T_0)f\right\| + \left\|(ST_0 - S_0T_0)f\right\|$$
$$\leq \left\|(T - T_0)f\right\| + \left\|(S - S_0)T_0f\right\| < \frac{1}{2}\varepsilon+\frac{1}{2}\varepsilon =\varepsilon,$$ 
i.e. $\left\{ST : (S,T)\in \mathcal W\right\} \subset \mathcal V$. 
To show that $T\mapsto T^{-1}$ is continuous in $T_0$ let 
$\mathcal V = \left\{S\in \mathrm U(\mathcal H): \left\|(S-T_0^{-1})f\right\| < \varepsilon\right\}$ a typical neighbourhood of $T_0^{-1}$ in $\mathrm U(\mathcal H)$. For $g := T_0^{-1}f$ let $T\in \mathrm U(\mathcal H)$ satisfy $\left\|(T-T_0)g\right\| < \varepsilon$. Then
$$\left\|(T^{-1}-T_0^{-1})f\right\| = \left\|T^{-1}T_0g - g\right\| = \left\|T_0g - Tg\right\| < \varepsilon,$$ 
i.e. $\left\{T^{-1} : \left\|(T-T_0)g\right\| < \varepsilon\right\} \subset \mathcal V$. \qed\\

This result with its simple proof is only worthwhile to publish because in the literature at several places the contrary is stated and because therefore some extra but superfluous efforts have been made. For example, Simms \cite{Si} explicitly states that the unitary group is not a topological group in its strong topology and that therefore the proof of Bargmann's theorem \cite{Ba} has to be rather involved. But also recently in the paper of Atiyah and Segal \cite{AS} some proofs and considerations are overly complicated because they also state that the unitary group is not a topological group. The assertion of proposition 1 has been mentioned in \cite{Sc}.\\

The misunderstanding  that $\mathrm U(\mathcal H)$ is not a topological group in the strong topology might come from the fact that the composition map $$\mathrm B(\mathcal H) \times \mathrm B(\mathcal H) \to \mathrm B(\mathcal H), (S,T) \mapsto ST,$$ is not continuous in the strong topology and consequently  $\mathrm{GL} (\mathcal H)$ is not a topological group (in the infinite dimensional case). But the restriction of the composition to  $\mathrm U(\mathcal H) \times \mathrm U(\mathcal H)$ is continuous since all subsets of $\mathrm U(\mathcal H)$ are uniformly bounded and equicontinuous.\\

Another assertation in \cite{AS} is that the compact open topology on $\mathrm U(\mathcal H)$ does not agree with the strong topology and therefore some efforts are made in \cite{AS} to overcome this assumed difficulty. However, again because of the uniform boundedness of the operators in $\mathrm U(\mathcal H)$ one can show:\\

\textbf{Proposition 2:} \textit{The compact open topology on $\mathrm U(\mathcal H)$ coincides with the strong topology.}\\

\textbf {Proof.} The compact open topology on $\mathrm{B}(\mathcal H)$ and hence on $\mathrm U(\mathcal H)$ is generated by the seminorms $T \mapsto \left\|T\right\|_K := \mathrm{sup}\left\{\left\|Tf\right\| : f \in K\right\}$ where $K \subset \mathcal H$ is compact. Let $\mathcal V = \left\{T \in \mathrm U(\mathcal H) : \left\|T - T_0\right\|_K < \varepsilon \right\}$ be a typical neighbourhood of $T_0\in \mathrm U(\mathcal H)$ where $K \subset \mathcal H$ is compact and $\varepsilon > 0$. We have to find a strong neighbourhood $\mathcal W$ of $T_0$ such that $\mathcal W \subset \mathcal V$. Let $\delta := \frac 1 3 \varepsilon$. By compactness of $K$ there is a finite subset $F\subset \mathcal H$ such that $K \subset \bigcup\left\{B(f,\delta) : f \in F\right\}$ where $B(f,r) = \left\{g\in \mathcal H : \left\|f - g\right\| < r \right\}$ is the usual open ball around $f$ of radius $r$. Now, for $k \in K$ there exist $f\in F$ with $k\in B(f,\delta)$ and $g\in B(0,\delta)$ such that $k = f + g$. We conclude, for $\left\|T-T_0\right\|_F < \delta$
$$\left\|(T-T_0)k\right\| \leq \left\|(T-T_0)f\right\| + \left\|(T-T_0)g\right\| < \delta + 2\delta = \varepsilon.$$
As a consequence, the strongly open $\mathcal W = \left\{T \in \mathrm U(\mathcal H) : \left\|T-T_0\right\|_F < \delta \right\}$ ist contained in $\mathcal V$. \qed\\

This proof essentially shows that on an equicontinuous subset of $\mathrm B(\mathcal H)$ the strong topology is the same as the compact open topology and both topologies coincide with the topology of uniform convergence on a total subset $D\subset \mathcal H$. \\

In particular, if $\mathcal H$ is separable with orthonormal basis $(e_k)_{k\in \mathbb N}$, the seminorms $T \mapsto \left\|Te_k\right\|$ generate the strong topology. A direct consequence is (in contrast to an assertion in Wikipedia\footnote{\tt\scriptsize{http://en.wikipedia.org/wiki/Metrization\_theorem\#Examples\_of\_non-metrizable\_spaces} (30.08.13)} which explicitly presents $\mathrm U(\mathcal H)$ with respect to the strong topology as an example of a non-metrizable space): \\

\textbf{Proposition 3:} \textit{The strong topology on $\mathrm U(\mathcal H)$ is metrizable and complete if $\mathcal H$ is separable.}\\

$\mathrm U(\mathcal H)$ is complete since the limit of a sequence of unitary operators which converges pointwise is again unitary.\\

The remarkable result of Kuiper \cite{Ku} that $\mathrm U(\mathcal H)$ is contractible in the norm topology if $\mathcal H$ is infinite dimensional and separable has been generalized by \cite{AS} to $\mathrm U(\mathcal H)$ with the compact open topology. By  proposition 2 we thus have\\

\textbf{Corollary:} \textit{$\mathrm U(\mathcal H)$ is contractible in the strong topology if $\mathcal H$ is infinite dimensional and separable.}\\

\textbf{Remark} \textit{The first three results extend to the projective unitary group $\mathrm{PU}(\mathcal H)= \mathrm U(\mathcal H)/\mathrm U(1) \cong \mathrm U(\mathbb P\mathcal H)$: This group is again a topological group in the strong topology, the strong topology coincides with the compact open topology and it is metrizable and complete for separable $\mathcal H$. Moreover we have the following exact sequence of topological groups 
$$1 \longrightarrow\mathrm U(1) \longrightarrow \mathrm U(\mathcal H) \longrightarrow \mathrm{PU}(\mathcal H) \longrightarrow 1$$ exhibiting $\mathrm{U}(\mathcal H)$ as a central extension of $\mathrm{PU}(\mathcal H)$ by $\mathrm U(1)$ in the context of topological groups and at the same time as a $\mathrm U(1)$-bundle over $\mathrm{PU}(\mathcal H)$.}\\

As a consequence, $\mathrm{PU}(\mathcal H)$ is simply connected (with respect to the strong topology), but not contractible. $\mathrm{PU}(\mathcal H)$ is an Eilenberg-MacLane space $K(\mathbb Z,2)$.\\ 

The sequence is not split as an exact sequence of topological groups or as an exact sequence of groups. Moreover, one can show that there does not exist a continuous section $\mathrm{PU}(\mathcal H) \to \mathrm{U}(\mathcal H)$ \cite{Sc}.\\

In view of the result of proposition 1 it is natural to ask whether $\mathrm{U}(\mathcal H)$ has the structure of a Lie group with respect to the strong topology. Let us review what happens in the case of the norm topology:\\

We know that $\mathrm{U}(\mathcal H)$ is a real Banach Lie group in the norm topology: Its local models are open subsets of the space $L \subset \mathrm B(\mathcal H)$ of bounded skew-symmetric operators. $L$ is a real Banach space and a Lie algebra with respect to the commutator. The exponential map
$$\exp : L \to \mathrm{U}(\mathcal H), B\mapsto \exp B = \sum \frac {B^n} {n!}$$
is locally invertible and thus provides the charts to define the manifold structure on the unitary group. In this way, 
$\mathrm{U}(\mathcal H)$ is a Lie group with Lie algebra $L$.\\

The same procedure does not work for the strong topology except for $\dim \mathcal H < \infty$. Although it can be shown that the above exponential map $\exp: L \to \mathrm{U}(\mathcal H)$ is continuous with respect to the strong topologies, it is not a local homeomorphism. Another way to see that $\mathrm{U}(\mathcal H)$ cannot be a Lie group with local models in $L$ with respect to the strong topology was told to me by K.-H. Neeb: Choose an orthonormal basis in $\mathcal H$. The subgroup $K \subset \mathrm{U}(\mathcal H)$ of diagonal operators with respect to an orthonormal basis $(e_j)_{j\in J}$ is the abelian group $$K = \left\{T \in \mathrm{U}(\mathcal H) : T = (\lambda_j)_{j\in J}, \left|\lambda_j\right|=1\right\}$$ isomorphic to the product of infinitely many circles $\mathrm U(1)$. The topology on $K$ induced from the strong topology is the product topology. Hence, $K$ ist compact. If $\mathrm{U}(\mathcal H)$ would be a Lie group in the strong topology then $K$ would be a Lie group as well with models in the space $D$ of diagonal operators in $L$ (with the strong topology). However, as a compact Lie group $K$ would have to be a finite dimensional manifold.\\

Note that if $\exp$ were locally invertible for the strong topologies then the same would be true for the restriction $\exp: D \to K$. But this restriction is not locally invertible, since for every neighbourhood $\mathcal V \subset K$ of $1=\mathrm{id}_{\mathcal H}$ the inverse image $\exp^{-1}(\mathcal V)$  contains all but finitely many straight lines of the form 
$$\mathbb R_{m} := \left\{T\in D : \,T = (\lambda_j), \,\lambda_j = 0 \,\,\, \mathrm{for} \,\,\,j\neq m, \,\lambda_m\in i\mathbb R \right\} \cong \mathbb R,$$ 
where $m\in \mathbb N$, and $\exp$ is not injective on $\mathbb R_m$. \\
 
According to  its importance in mathematics and physics one might be tempted to use all the unitary, strongly continuous  one parameter groups $$\mathbb R \to \mathrm U(\mathcal H), t \mapsto B(t), t\in \mathbb R,$$  as the basic geometric and analytic information to find a manifold structure on $\mathrm{U}(\mathcal H)$. By Stone's theorem these are exactly the one parameter groups $$t \mapsto \exp itA, t\in \mathbb R,$$ for self adjoint (not necessarily bounded) operaters $A$ on $ \mathcal H$. However, the set of all self adjoint operators is not a linear space.\\

%\texttt{variant 1}

%We finally explain why it is important to know that our groups $\mathrm{U}(\mathcal H)$ and $\mathrm{PU}(\mathcal H)$ is a topological group for example in the context of studying Hilbert bundles. Roughly, this is the case because we can construct the classifying spaces $B\mathrm{U}(\mathcal H)$ and $B\mathrm{PU}(\mathcal H)$ using proposition 1. \\

%\texttt{variant 2}

The result of proposition 1 that $\mathrm U(\mathcal H)$ with the strong topology is a topological group helps to find simpler and more transparent proofs (e.g. than those in \cite{Si} and \cite {AS}) and it gives a coherent picture when dealing with topological fiber bundles or with unitary representations of topological groups. 
%In particular, it is of structural significance in the study of Hilbert bundles, since the associated frame bundles obtain the structure of a principal fiber bundle with structural group $\mathrm U(\mathcal H)$ or $\mathrm {PU}(\mathcal H)$ and it allows to construct the related classifying spaces $B\mathrm{U}(\mathcal H)$ and $B\mathrm{PU}(\mathcal H)$ using proposition 1. Let us explain this latter situation in some detail:\\
%\texttt{variant 3}
In the following we exemplify the advantage of knowing that $\mathrm U(\mathcal H)$ is a topological group with respect to the strong topology by applying this result to the study of Hilbert bundles. For a given topological group $G$ the homotopy classification of all equivalence classes of principal fiber bundles over a fixed paracompact space $X$ can be described using the classifying space $BG$. The significance of proposition 1 is that this can also be done for $G=\mathrm U(\mathcal H)$ or $\mathrm{PU}(\mathcal H)$ with the strong topology. Let us explain the consequences for the study of Hilbert bundles:  \\

A Hilbert bundle $E$ over a (paracompact) space $X$ is a bundle $\pi: E\to X$ over $X$ with continuous projection $\pi$ such that the fibers $E_x = \pi^{-1}(x), x\in X,$ are isomorphic to a separable complex Hilbert space $\mathcal H$ or its projectivation $\mathbb P(\mathcal H)$. Here, 'isomorphic' means unitarily isomorphic. In particular, there exists a cover of open subsets $V \subset X$ with bundle charts (i.e. homeomorphisms)
$$\phi: E|_V \to V \times \mathcal H$$
(when of $\mathcal H \cong E_x$ as the typical fiber) such that 
$\mathrm{pr}_1 \circ \phi = \pi$ and 
$$\phi_x := \mathrm{pr}_2\circ \phi|_{E_x} : E_x \to \mathcal H$$ is unitary for all $x\in X$. 
Thus, for $\dim \mathcal H = n < \infty$ the bundle $E$ is an ordinary complex vector bundle with typical fiber $\mathbb C^n$ and structural group $\mathrm U(n)$.\\ 

The transition maps for further bundle charts 
$\phi': E|_V' \to V' \times \mathcal H$, $W=V\cap V' \neq \emptyset$, are
$$ \phi'\circ \phi^{-1}: W\times \mathcal H \longrightarrow W\times \mathcal H$$
which are completely determined by their projections
$$ \psi = \psi_{V',V} :=\mathrm{pr}_2\circ \phi'\circ \phi^{-1}: W\times \mathcal H \longrightarrow  \mathcal H.$$
Now, as we have shown above, $\psi$ is continuous, if and only if the induced map
$$\widehat \psi: W \to \mathrm U(\mathcal H), x \mapsto (f \mapsto \psi(x,f)),$$
is strongly continuous.  $\widehat \psi$ will not be continuous with respect the norm topology, in general. An example will be given below. In the case of $\mathbb P(\mathcal H)$ as the typical fiber of $E$ we have analoguous statements.\\

As a consequence, the natural principal fiber bundle $P = P_E \to X$ associated to the Hilbert bundle $E$ (the frame bundle with fibers $P_x = \mathrm U(E_x,\mathcal H)$ if $\mathcal H$ is the typical fiber) will be a principal fiber bundle whose structural group is $\mathrm U(\mathcal H)$ with its strong topology and, in general, not with respect to the norm topology. Note that $P_E$ will be, in addition, a principal fiber bundle with respect to the norm topology on $\mathrm{U}(\mathcal H)$ if and only if there exists a cover of $X$ with bundle charts such that all the induced transition maps $\widehat \psi : W \to \mathrm{U}(\mathcal H)$ are norm continuous. Let us call such a bundle 'norm-defined'. \\

In the case that $\mathbb P (\mathcal H)$ is the typical fiber of $E$ we have the analoguous results for the associated principal bundle $P_E$ (with fibers $P_x = \mathrm U(E_x, \mathbb P\mathcal H)$): The structural group is $\mathrm{PU}(\mathcal H)$ with the strong topology in general. Morever, whenever $E$ is norm-defined $P_E$ can also be viewed as to be a principal fiber bundle with structural group the projective unitary group $\mathrm{PU}(\mathcal H)$ in its norm topology.\\

In order to classify the Hilbert bundles over $X$ it is enough to classify the principal fiber bundles with structural groups $\mathrm U(\mathcal H)$ resp. $\mathrm {PU}(\mathcal H)$. Let $Princ_{\mathrm U(\mathcal H)}^N(X)$ the set of isomorphism classes of principal fiber bundles with $\mathrm U(\mathcal H)$ in the norm topology and correspondingly $Princ_{\mathrm U(\mathcal H)}^S(X)$ the set of isomorphism classes of principal fiber bundles with $\mathrm U(\mathcal H)$ in the strong topology. Analoguously, we define $Princ_{\mathrm {PU}(\mathcal H)}^A(X)$ for $A \in \left\{N, S\right\}$.\\

Unitary group (vector bundles): Since the unitary group is contractible in both topologies every principal bundle is trivial: $$Princ_{\mathrm U(\mathcal H)}^A(X) \cong [X, B\mathrm U(\mathcal H)^A] = \left\{[X\times\mathcal \mathrm U(\mathcal H)]\right\}$$ for  $A \in \left\{N, S\right\}$. For an arbitrary Hilbert bundle with typical fiber $\mathcal H$ this implies that it is already isomorphic to the trivial bundle $X\times \mathcal H$. For the norm-defined bundles $E$ the associated principal bundle $P_E$ is in $Princ_{\mathrm U(\mathcal H)}^N(X)$ and an isomorphism $E \cong X\times \mathcal H$ can be found which is locally given by transition functions which are induced by norm continuous $W\to \mathrm U(\mathcal H)$. Note, that the classifying spaces $B\mathrm{U}(\mathcal H)^{A}$ are weakly contractible (all homotopy groups are trivial) for $A\in \left\{N, S\right\}$.\\
 
Projective unitary group (projective bundles): We know already that $\mathrm{PU}(\mathcal H)$ is a $K(\mathbb Z, 2)$ for both topologies on the projective unitary group which we will indicate by a superscript $A$. From the homotopy sequence  corresponding to the universal bundle
$$\mathrm{PU}(\mathcal H)^{A} \longrightarrow \mathrm{EPU}(\mathcal H)^{A} \longrightarrow \mathrm{BPU}(\mathcal H)^{A}$$ 
one concludes that $\mathrm {BPU}(\mathcal H)^{A}$ is an Eilenberg-MacLane space $K(\mathbb Z, 3)$. Now, the homotopy classification of principal fiber bundles asserts that there is a bijection between  $Princ_{\mathrm {PU}(\mathcal H)}^A(X)$ and $[X,\mathrm{BPU}(\mathcal H)^{A}]$, the set of homotopy classes of continuous $X\to \mathrm{BPU}(\mathcal H)^{A}$. For a $K(\mathbb Z,3)$ this is cohomology:  $[X,\mathrm{BPU}(\mathcal H)^{A}] \cong H^3(X,\mathbb Z)$. We arrive at the following result which is essentially contained in a different form in \cite{SP}: \\

\textbf{Proposition 4:} \begin{itemize}
\item \textit{The isomorphism classes of  Hilbert bundles over X are in one-to-one correspondence to $H^3(X,\mathbb Z) \cong [X,B\mathrm{PU}(\mathcal H)^S]\cong Princ_{\mathrm{PU}(\mathcal H)}^S(X)$.}
\item \textit{The isomorphism classes of norm-defined  Hilbert bundles over X are also in one-to-one correspondence to $H^3(X,\mathbb Z) \cong [X,B\mathrm{PU}(\mathcal H)^N]\cong Princ_{\mathrm{PU}(\mathcal H)}^N(X)$ where the isomorphisms of the Hilbert bundles are given by norm continuous transition maps.}
\end{itemize}

Note, that the zero element of $H^3(X,\mathbb Z)$ represents the class of all trivial bundles which also can be described as the projective Hilbert bundles $E$ of the form $\mathbb P F$ where $F$ is a true vector bundle with fibers $F_x \cong \mathcal H$. Proposition 4 also implies that in every equivalence class of Hilbert bundles there exists a norm-defined representative.\\

\textbf{Example:} Let $\pi: P\to X$ be a principal fiber bundle with structural group a compact group $G$ and let $\mu$ be the normalized Haar measure. For $\mathcal H = L_2(G, \mu)$ let $L : G\times \mathcal H \to \mathcal H $ the induced left action (regular representation). The example of a possibly non norm-defined Hilbert bundle is the associated $E:= P \times_G \mathcal H$ with structural group the unitary group with the strong topology. The bundle charts of $E$ coming from bundle charts of $P$ and describung the structure of $Q$ completely are the following: If $\phi: P|_V\to V\times G$ is a bundle chart, then $\psi: Q|_V \to V\times \mathcal H$ is given by $[(p,f)] \mapsto (\pi(p), L_{\mathrm{pr}_2\circ \phi(p)}(f)$. As a consequence, the transition maps for these charts cannot be induced be norm continuous maps $W \to \mathrm U(\mathcal H)$, in general.\\

Finally, let us remark, that one reason to search for a Lie group structure on $\mathrm U(\mathcal H)$ with the strong topology is that one wants to study smooth Hilbert bundles over (finite dimensional) manifolds $X$: These are Hilbert bundles in the above sense and smooth manifolds with models in a Hilbert space where the transition maps $W\times \mathcal H \to \mathcal H$ or  $W\times \mathbb P \mathcal H \to \mathbb P \mathcal H$ ($W\subset X$ open) are smooth. Of course, one expects the associated frame bundle to be smooth as well with respect to the not yet defined Lie group structure on $\mathrm U(\mathcal H)$ or $\mathrm {PU}(\mathcal H)$. In particular, we have the natural question in what sense the corresponding strongly continuous $W\to \mathrm U(\mathcal H)$ or $W\to \mathrm{PU}(\mathcal H)$ are differentiable.

 \vspace{0,5cm}

\end{document}